\documentclass[12pt]{article}
\usepackage{graphicx} 
\usepackage{amsmath} 
\usepackage{amssymb} 
\usepackage{amsfonts} 
\usepackage{amsthm} 
\usepackage{tikz}

\usepackage{color}
\usepackage{subfig}
\input{epsf}
\input xy
\xyoption{all}
\usepackage{latexsym}
\newtheorem{lemma}{Lemma}[section] 
\newtheorem{theorem}[lemma]{Theorem} 
 
\newtheorem{corollary}[lemma]{Corollary}

\theoremstyle{definition}
\newtheorem{definition}[lemma]{Definition}
\newtheorem{example}[lemma]{Example}
\newtheorem{remark}[lemma]{Remark}

\renewcommand{\epsilon}{\varepsilon} 

\renewcommand{\leq}{\leqslant} 
\renewcommand{\geq}{\geqslant}

\newcommand{\Po}{\mathcal{P}}
\newcommand{\R}{\mathcal{R}}
\newcommand{\Q}{\mathcal{Q}} 
\newcommand{\F}{\mathcal{F}}

\newcommand{\Mon}{\mathrm{Mon}} 

\newcommand{\St}{\mathrm{Stab}}

\newcommand{\K}{\mathcal{K}}

\newcommand{\Ga}{\Gamma}

\newcommand{\rk}{\mathrm{rank}} 
\begin{document}


\title{Finite Polytopes have Finite Regular Covers}
 
\author{B. Monson\thanks{Supported in part by PAPIIT-Mexico grant \# IN112512
and by  NSERC of Canada Discovery Grant \# 4818
}\\
University of New Brunswick\\
Fredericton, New Brunswick, Canada E3B 5A3
\and and\\[.05in]
Egon Schulte\thanks{Supported by NSF-grant DMS--0856675}\\
Northeastern University\\
Boston, Massachussetts,  USA, 02115}

\date{ \today }
\maketitle

\begin{abstract}
\noindent
We prove that any finite, abstract $n$-polytope is covered by a finite,
abstract regular  $n$-polytope.

\medskip


\end{abstract}

\section{Introduction}
Any point or line segment is a regular convex polytope, admittedly of 
modest dimension. And every polygon with $p$ vertices
is combinatorially equivalent to a regular $p$-gon. Moving 
further into the combinatorial domain, it is clear that
a typical map $\mathcal{M}$
on a compact surface is  non-regular. Even so, it is folklore (and
not so hard to prove) that
$\mathcal{M}$ is a quotient of a  regular map $\R$,  
most likely on some new surface,
and which from another point of view covers the base map $\mathcal{M}$.

On the other hand, it is certainly true that any  $n$-polytope 
is a quotient of the universal $n$-polytope 
$\mathcal{U}_n = \{\infty, \ldots,\infty\}$,
which is indeed infinite if the rank $n\geq 2$. 

In this note we look at all these situations from an abstract point of view
and prove  a reasonable but  hitherto unaddressed result:
every finite abstract $n$-polytope $\Q$ has a finite regular cover $\R$
(Theorem~\ref{main}).

\section{Abstract Polytopes and their Automorphism Groups}\label{abst}

An abstract  $n$-polytope $\Po$  has certain key 
combinatorial properties of the face lattice of a convex $n$-polytope;
in general, however, $\Po$ need not be a lattice or have
a familiar geometric realization.
Let us  summarize some general definitions and results, 
referring  to \cite{arp} for details.
An \textit{abstract $n$-polytope} $\Po$ is a partially ordered
set with  properties \textbf{A}, \textbf{B} and 
\textbf{C} below.

\medskip

\noindent\textbf{A}: $\Po$ has a strictly monotone rank function with range 
$\{-1,0,\ldots,n\}$. Moreover,
$\Po$ has a unique least face
$F_{-1}$ and   unique greatest face $F_n$. 

\medskip\noindent
An element $F \in \Po$ with $\rk(F)=j$ is
called a $j$-\textit{face}; often  $F_j$ will indicate a $j$-face. 
Each maximal chain or \textit{flag} in $\Po$ therefore contains $n+2$ faces,
so that $n$ is the number of \textit{proper} faces in each flag.
We let $\F(\Po)$ be the set of all flags in $\Po$.
Naturally, faces of ranks 0, 1 and $n-1$ are called 
vertices, edges and  facets, respectively.
\medskip 

\noindent\textbf{B}: Whenever $F < G$ with $\rk(F)=j-1$ and
$\rk(G)=j+1$, there are exactly two $j$-faces $H$ with
$F<H<G$.

\medskip\noindent
For $0 \leq j \leq n-1$ and any flag $\Phi$, 
there thus exists a unique \textit{adjacent} 
flag $\Phi^j$, differing from $\Phi$ in just the
face of rank $j$ . With this notion of adjacency,  $\F(\Po)$
becomes the \textit{flag graph} for $\Po$.
If $F \leq G$ are incident faces  
in  $\Po$,  we call
$$ G/F := \{ H \in \mathcal{P}\, | \, F \leq H \leq G \}\;$$
a \textit{section}  of $\Po$.
 
\medskip
\noindent\textbf{C}: $\Po$ is \textit{strongly flag--connected}, 
that is, the flag graph for each section is connected. 
 
\medskip\noindent
It follows that 
$G/F$ is a  ($k-j-1$)-polytope in its own right, if 
$F \leq G$ with $\mathrm{rank}(F) = j \leq k = \mathrm{rank}(G)$.
For example, if  $F$ is a vertex, then the
section $F_n/F$ is called the \textit{vertex-figure} over  $F$.
Likewise, it is useful to think of the  $k$-face $G$ 
as having the structure of the $k$-polytope $G/F_{-1}$.

The \textit{automorphism group} $\Gamma(\Po)$ consists of all
order-preserving bijections on $\Po$. We say  $\Po$ is \textit{regular} if 
$\Gamma(\Po)$ is transitive on the flag set $\F(\Po)$. In this case we
may choose any one flag $\Phi \in \F(\Po)$ as \textit{base flag},  
then  define $\rho_j $ to be  the (unique) automorphism  
mapping $\Phi$ to $\Phi^j$, for $0 \leq j \leq n-1$. Each $\rho_j$ has period $2$.
 From \cite[2B]{arp}
we recall that $\Gamma(\Po)$ is then a \textit{string C-group}, meaning that 
it  has  the following properties \textbf{SC1} and \textbf{SC2}:

\medskip\noindent\textbf{SC1}:   $\Gamma(\Po)$ is 
generated by $ \{\rho_0, \ldots, \rho_{n-1}\}$. These  involutory generators
satisfy the commutativity relations
typical of a Coxeter group with string diagram, namely  
\begin{equation}\label{relreg}
(\rho_ j \rho_k)^{p_{jk}} = 1,\;  \mathrm{ for }\; 
0 \leq j \leq k \leq n-1, 
\end{equation}
where   $p_{jj} = 1$ and $p_{jk} = 2$ whenever $|j-k|>1$.
In other words, $\Gamma(\Po)$ is a
\textit{string  group generated by involutions} 
or \textit{sggi}.
 
\medskip\noindent\textbf{SC2}: $\Gamma(\Po)$ satisfies the
\textit{intersection condition}
\begin{equation}\label{interreg}
\langle \rho_k : k \in  I\rangle  \cap \langle \rho_k : k \in  J\rangle  = 
\langle \rho_k : k \in  I \cap J\rangle, 
\;  \mathrm{ for\,\, any }\; I, J \subseteq \{0 \ldots, n-1\}\;.
\end{equation}

\medskip

\noindent
The fact that one can reconstruct 
a regular polytope in a canonical way 
from any string C-group $\Gamma$ is at the heart of the theory \cite[2E]{arp}.

The periods $p_j := p_{j-1,j}$ in (\ref{relreg}) 
satisfy $2 \leq p_j \leq \infty$ and        
are assembled into
the \textit{Schl\"{a}fli  symbol}  $\{ p_1, \ldots, p_{n-1}\}$ for the
regular polytope $\Po$. We also say that  $\Po$  has  \textit{type} indicated
by the same symbol and that the group $\Ga(\Po)$ has rank $n$.

In the same way, any sggi $G  = \langle g_0, \ldots, g_{n-1} \rangle$
also has a Schl\"{a}fli symbol, type and rank, 
although $p_j = 1$ can occur (if $g_{j-1} = g_j$; this is impossible for regular polytopes
by (\ref{interreg})).

The \textit{dual} $\Po^\ast$ of the polytope $\Po$ is obtained by simply reversing
the partial order on the underlying set of faces. If $\Po$ is regular of type  
$\{ p_1, \ldots, p_{n-1}\}$, then $\Po^\ast$ is also regular of type
$\{ p_{n-1}, \ldots, p_{1}\}$.

\begin{definition} \cite[2D]{arp}   
Let $\R$ and $\Po$ be $n$-polytopes. A \textit{covering}
is a rank and adjacency preserving homomorphism 
$\eta : \R \rightarrow \Po$. 
(This means that $\eta$ induces a mapping 
$\F(\R) \rightarrow \F(\Po)$
which sends any $j$-adjacent pair of flags in $\R$ 
to another such pair in $\Po$; it is easy to show that $\eta$ must then be surjective.) 
We also say that
$\R$ is a \textit{cover} of $\Po$ and write $\R \rightarrow \Po$.
\end{definition}

\medskip

If $\R$ covers $\Po$, then from another point of view, $\Po$
will be a quotient of $\R$; see \cite[2D]{arp}, \cite{hartley5} or \cite{mixA}.
One way  to understand how $\Po$ arises  by identifications in $\R$ 
is to exploit the monodromy group:

\begin{definition} Let $\Po$ be a polytope of rank $n \geq 1$. For
$0 \leq j \leq n-1$, let $g_j$ be the bijection on $\F(\Po)$
which maps each flag $\Phi$ to the $j$-adjacent flag $\Phi^j$. Then the 
\textit{monodromy group} for $\Po$ is
$$ \Mon(\Po) = \langle g_0, \ldots, g_{n-1}\rangle\;$$
(a subgroup of the symmetric group on $\F(\Po)$).
\end{definition}

It is easy to see that $\Mon(\Po)$ is an sggi.
Let us quote from \cite{mixA} some useful and fairly easily proved results.
We let $\rm{Stab}_{\rm{Mon}(\mathcal{P})}\Phi$ denote the stabilizer of a 
flag $\Phi$ under the action of $\rm{Mon}(\mathcal{P})$.

 \begin{theorem}\label{monregisomorphism}
 
Let $\R$ be a regular $n$-polytope with  base flag $\Phi$,
automorphism group $\Gamma(\R) = \langle  \rho_0, \ldots , \rho_{n-1}\rangle$,
and monodromy group $\Mon(\R) = \langle g_0, \dots, g_{n-1} \rangle$. Then 
there is an isomorphism  $\Gamma(\R) \simeq   \Mon(\R) $
mapping each $\rho_j$ to $g_j$.
\end{theorem} 

 \medskip

\newpage

\begin{theorem}\label{epistocover}
Suppose that $\R$ and $\Q$ are $n$-polytopes and that
$$\bar{\eta} : \Mon(\R) \rightarrow \Mon(\Q)$$
is an epimorphism of sggi's 
(i.e. mapping specified generators to specified generators, in order).
Suppose also that there are flags 
$\Lambda'$ of $\R$ and $\Lambda$ of $\Q$ such that
\begin{equation}\label{twostabs}
(\St_{\Mon(\R)}\Lambda')\, \bar{\eta} \subseteq \St_{\Mon(\Q)}\Lambda\; .
\end{equation}
Then there is a unique covering $\eta : \R \rightarrow \Q$
which maps $\Lambda'$ to $\Lambda$.
\end{theorem}

\medskip

\begin{remark}
If $\R$ is regular, then condition 
(\ref{twostabs}) is fulfilled automatically, since all flags 
$\Lambda'$ are equivalent, with trivial stabilizer, in 
$\Gamma(\R) \simeq \Mon(\R)$ (Theorem~\ref{monregisomorphism}).
\end{remark}

\section{Regular covers of general polytopes}

A key step in our construction is provided by the following

\begin{lemma}\label{2K}
 Suppose $\K$ is a regular $m$-polytope with automorphism group 
$\Ga(\K) = \langle \rho_0, \ldots, \rho_{m-1} \rangle$, 
Schl\"{a}fli type $\{p_1. \ldots, p_{m-1}\}$ and facet set $A$. 
Then there is a regular $(m+1)$-polytope $\bar{\K}$ of type 
$\{p_1, \ldots, p_{m-1},4\}$, with facets isomorphic to
$\K$, and such that 
$|\Ga(\bar{\K})| = |\Ga(\K)|\cdot 2^{|A|}$, when $\mathcal{K}$ is finite.
Indeed, $\bar{\K}$ is finite if and only if $\K$ is finite.
\end{lemma}

\medskip
\noindent\textbf{Proof}. We could simply refer to \cite[Theorem 8C2]{arp} and take
$\bar{\K} = {(2^{\K^\ast})}^{\ast}$. But for the reader's convenience, we shall 
rework here the essentials of that construction.
First off, let $N$ be the group of all sequences
$$x: A \rightarrow C_2 := \{\pm 1\}, $$
with component-wise multiplication. 
(If $A$ is infinite, take sequences of finite ``support'',
meaning that $-1$ can occur only finitely often.)
Thus $N$ is generated by the indicator functions $\delta_F, F \in A$, where
for facets $F, H \in A$ we have
$$ (H)\delta_F = \left\{
\begin{array}{ll}
 -1,& \mathrm{if}\; H = F,\\
 +1,& \mathrm{if}\; H \neq F.
\end{array}
\right. $$
We require $\rho_m := \delta_{F_{m-1}}$, the indicator function for 
the base facet of $\K$.

Now for $x\in N, \gamma \in \Ga(\K), F \in A$
we let 
$$(F) x^{\gamma} := (F)\gamma^{-1} x.$$
This defines an action $x \mapsto x^{\gamma}$ of $\Ga(\K)$ on $N$, and we let
$$ S:= N \rtimes \Ga(\K)$$
be the corresponding semidirect (indeed, wreath) product. For ready computation
 we abuse notation a bit, taking $N\triangleleft S$ and $\Ga(\K) < S$, so that
$x\cdot\gamma = \gamma\cdot x^{\gamma}$, for $x \in N, \gamma \in \Ga(\K)$.
Note in particular that $\delta_F^{\gamma} = \delta_{(F)\gamma}$, for any facet $F$.
Since the base facet $F_{m-1}$
is stabilized by $\rho_j, 0\leq j \leq m-2$, we have 
$$\rho_j \rho_m = \rho_j \delta_{F_{m-1}}
= \delta_{F_{m-1}}^{\rho_j}\rho_j = \delta_{(F_{m-1})\rho_j} \rho_j
= \delta_{(F_{m-1})} \rho_j = \rho_m \rho_j\;.
 $$
Similarly we find that $(\rho_{m-1}\rho_m)^2$ is an element of $N$ and so has period $2$.
In fact,  $$ S = \langle \rho_0, \ldots, \rho_{m-1},\rho_m \rangle$$
is an sggi of rank $m+1$ and type $\{p_1, \ldots, p_{m-1},4\}$.
Using the unique factorization in $S$ given by $S = N \,\Ga(\K)$, together with the 
intersection property (\ref{interreg}) for $\K$, we soon verify the intersection 
property for $S$, too \cite[Lemma 8B5]{arp}. Thus $S$ is a string C-group of rank $m+1$, so
$S = \Ga(\bar{\K})$ for just the sort of regular $(m+1)$-polytope $\bar{\K}$ which we seek.
\hfill $\square$
\medskip

\begin{remark}
In  \cite{pellicer2},  Daniel Pellicer uses   `CPR-graphs' to 
generalize the results in Lemma~\ref{2K}. One can construct
a  regular polytope  
$\bar{\K}$ of type 
$\{p_1, \ldots, p_{m-1},2s\}$, with facets isomorphic to
$\K$, for any   integer $s\geq 2$.
\end{remark}

\medskip
\begin{lemma}\label{bump}
Let $G=\langle g_0,\ldots,g_{n-1}\rangle$ be an sggi, and let $0\leq i\leq n-2$. 
Suppose there exists a string C-group $H=\langle h_0,\ldots,h_{i}\rangle$ which 
covers $\langle g_0,\ldots,g_{i}\rangle$.  
Then there also exists a string C-group $L=\langle l_0,\ldots,l_{i+1}\rangle$ 
which covers $\langle g_0,\ldots,g_{i+1}\rangle$. Furthermore, we may choose $L$ in 
such a way that, if $j\leq i$ and 
$\langle g_0,\ldots,g_{j}\rangle\simeq\langle h_0,\ldots,h_{j}\rangle$, then 
$\langle g_0,\ldots,g_{j}\rangle\simeq\langle l_0,\ldots,l_{j}\rangle$ as well. 
Moreover, if $\langle g_0,\ldots,g_{i+1}\rangle$ and $H$ are finite, then we can 
take $L$ to be finite, too. 
\noindent
\end{lemma}
\medskip
\noindent\textbf{Proof}. Let $\K$ be a regular $(i+1)$-polytope with 
$\Ga(\K) \simeq H$. By Lemma~\ref{2K} there is a regular $(i+2)$-polytope $\bar{\K}$
with facets isomorphic to $\K$. We may suppose
$\Ga(\bar{\K}) =\langle h_0, \ldots, h_i, h_{i+1}\rangle$.

Now the \textit{mix} $L =  \langle g_0, \ldots, g_{i+1} \rangle \diamondsuit
\langle h_0, \ldots, h_{i+1} \rangle$ is the subgroup of the direct product 

\begin{equation}\label{dir}
\langle g_0, \ldots,  g_{i+1} \rangle \times
\langle h_0, \ldots,  h_{i+1} \rangle = \langle g_0, \ldots,  g_{i+1} \rangle \times
\Ga(\bar{\K})  
\end{equation}
generated by all $l_t := (g_t, h_t), 0\leq t \leq i+1$.
Clearly,  $L$ is  also an sggi of rank $i+2$.
(We refer to \cite{cunn1,arp, mixA,wilson2}
for other useful properties of this  operation.)

Since $H$ covers $\langle g_0, \ldots, g_{i} \rangle$,
we have $H \simeq \langle l_0, \ldots, l_{i} \rangle$, too.
Thus we can apply the quotient criterion \cite[Theorem 2E17]{arp}
to the second natural projection $L \rightarrow \Ga(\bar{\K})$ and conclude
that $L$ is a string C-group. 
The first natural projection shows that $L$ covers
$ \langle g_0, \ldots, g_{i+1}\rangle$. 

Appealing once more to Lemma~\ref{2K}, we see that
$L$ is finite if both factors in the direct product (\ref{dir})
are finite. 
   
\hfill $\square$
\medskip

The last lemma is just what we need to prove our main result.

\begin{theorem}\label{main}
{\rm (a)} Every finite sggi $G$ is covered by a finite string C-group $G'$.

{\rm (b)} Every finite $n$-polytope $\Q$  is covered by a finite 
regular $n$-polytope $\R$. If $\Q$ has all its $k$-faces isomorphic to 
some regular $k$-polytope  $\K$, then we may choose $\R$ 
to have its $k$-faces isomorphic to $\K$.
\end{theorem}

\noindent\textbf{Proof}. Let $G = \langle g_0, \ldots, g_{n-1} \rangle$ be any
finite sggi. Clearly, $H_0:=\langle g_0\rangle$ is a string C-group,
as is the dihedral group $H_1:=\langle g_0, g_1\rangle$. 
Thus we   can begin an 
inductive construction. Suppose $H_{j}$ is a finite string 
C-group covering $\langle g_0, \ldots, g_{j}\rangle$. Take $i=j$ and $H = H_{j}$
in    Lemma~\ref{bump}. We obtain a finite string C-group $L$ covering 
$\langle g_0, \ldots, g_{j+1}\rangle$. Now let $H_{j+1} = L$ and iterate.
Eventually we get a finite string C-group $G' = H_{n-1}$ which covers $G$.
Note  that if
$\langle g_0, \ldots, g_{k-1}\rangle$ happens to be a string C-group,   we 
can start the iteration with  this subgroup. It is clear from the Lemma~\ref{bump}
that we   end with the corresponding  subgroup of $G'$ unchanged up to 
isomorphism. 

For part (b) we merely apply part (a) to the finite sggi 
$G = \Mon(\Q) = \langle g_0, \ldots, g_{n-1}\rangle$.
Let $\R$ be the finite regular $n$-polytope  whose automorphism group is 
$G'$ constructed in (a). From Theorem~\ref{epistocover}
we conclude that $\R$ covers $\Q$. 
Since all $k$-faces are isomorphic to $\K$, it is also
true that $ \langle g_0, \ldots, g_{k-1}\rangle \simeq\Ga(\K)$ \cite{mixA}. 
We conclude that 
$\R$ has isomorphic $k$-faces.
\hfill $\square$

\begin{remark}
It is clear that a dual result concerning  co-$k$-faces must hold
in  Theorem~\ref{main}(b).
\end{remark}

\begin{corollary}
 Every convex $n$-polytope $\mathcal{Q}$ has a finite abstract regular 
cover $\mathcal{R}$. If $\mathcal{Q}$ is simplicial (or simple),
then $\mathcal{R}$ is likewise simplicial (or simple).
\end{corollary}

\medskip

The monodromy group  $G$ of a polytope $\mathcal{Q}$
always gives rise to a regular \textit{pre-polytopal\/} cover  
$\mathcal{V}$,
constructed from $G$ as a coset geometry in much the same way as a regular polytope
can be rebuilt from a given string C-group $\Gamma$. However, as the following example shows, 
this object $\mathcal{V}$ can fail condition \textbf{C} concerning strong flag-connectedness.

\begin{example}
Suppose $\Q$ is a pyramid  over the toroidal base
$\{4,4\}_{(3,0)}$. The `lateral' facets of this
  self-dual $4$-polytope  are the $9$ ordinary pyramids over 
the square faces in the toroid. Using GAP \cite{gap4}
we find that the monodromy group $G =  \langle g_0,g_1,g_2,g_3\rangle$
is an sggi of type $\{12,12,12\}$ and order $2^{12}\cdot3^{11}\cdot5$.
However, the intersection condition (\ref{interreg}) fails, since 
$ \langle g_1,g_2\rangle$ has index $2$ in 
$ \langle g_0,g_1,g_2\rangle\cap \langle g_1,g_2,g_3 \rangle$.
Following our earlier remarks, we could 
manufacture a regular \textit{pre-polytopal} cover
$\mathcal{V}$ of 
$\Q$ with automorphism group $G$. We find, however, that
the section between a typical vertex and facet of $\mathcal{V}$
consists of \textit{two} disjoint copies of a dodecagon $\{12\}$.

If we want a finite, regular \textit{polytopal} cover $\R$,
then we must appeal to Theorem~\ref{main}. Since the subgroup
$ \langle g_0,g_1,g_2\rangle$ is a string C-group (of order $2^{12}\cdot 3^3$),
we actually need to appeal to Lemma~\ref{2K} just once. The corresponding
regular $3$-polytope $\K$ has $4608$ facets $\{12\}$. Thus the regular extension
$\bar{\K}$ has type $\{12,12,4\}$ with group order 
$|\Ga(\bar{\K})| = 2^{12}\cdot 3^3\cdot 2^{4608}$.
The regular cover $\R$ still has type $\{12,12,12\}$, and its facets
are isomorphic to $\K$. The order of its automorphism group is bounded, somewhat absurdly,
by $$|G \times \Ga(\bar{\K})| = 2^{4632}\cdot 3^{14} \cdot 5\;. $$
Presumably a minimal regular cover of $\Q$ has much smaller group order.
\end{example}

\medskip
\begin{remark}
One might ask whether there is some sort of extension of Theorem~\ref{main}(b)
to the class of infinite polytopes $\mathcal{Q}$. Perhaps each has a regular cover with 
finite `covering index'.

To see that there is no hope for a general statement of this sort
 we begin with $\mathcal{T} = \{4,4\}$. This  familiar tiling of the Euclidean 
plane $\mathbb{R}^2$ by unit squares is an infinite regular $3$-polytope.
Next, for every odd integer $n\geq 3$, we subdivide the square with southwest vertex $(n,0)$,
using an $(n+2)$-gon together with $n$ triangles. 
The southeast 
vertex now has degree $n+3$. 
(The figure illustrates the case $n = 5$.)
The resulting $3$-polytope $\mathcal{Q}$ clearly has trivial
automorphism group $\Gamma(\mathcal{Q})$.

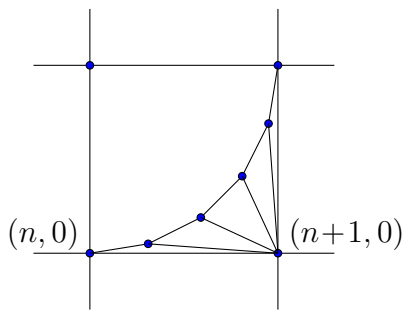
\begin{center}{\begin{tikzpicture}[scale=2.5]
\draw (-.3,0) -- (1.3,0);\draw (-.3,1) -- (1.3,1);
\draw (0,-.3) -- (0,1.3);\draw (1,-.3) -- (1,1.3);
\draw[fill=blue]  (0,0) circle (0.2mm);
\draw[fill=blue] (1,0) circle (0.2mm);
\draw[fill=blue] (0,1) circle (0.2mm);\draw[fill=blue] (1,1) circle (0.2mm);
\draw (0,0) -- (.31,.05) -- (.59, .19) -- (.81, .41) -- (.95, .69) -- (1,1);
\draw[fill=blue]  (.31,.05) circle (0.2mm);\draw[fill=blue] (.59, .19)  circle (0.2mm);
\draw[fill=blue]  (.81, .41) circle (0.2mm);\draw[fill=blue]  (.95, .69) circle (0.2mm);
\draw (1,0) -- (.31,.05);
\draw (1,0) -- (.59, .19);
\draw (1,0) -- (.81, .41);
\draw (1,0) -- (.95, .69);
\node [left] at (0.0,0.1) {$(n,0)$};\node [right] at (1.0,0.1) {$(n\!+\!1,0)$};
%
\end{tikzpicture}
}
\end{center}

Even though $\mathcal{Q}$ is locally finite, it is  also
 clear that any regular cover $\mathcal{R}$ must have 
Schl\"{a}fli type $\{\infty,\infty\}$.  But, in any case, each fibre
over a face of $\mathcal{Q}$ induced by the regular
cover $\mathcal{R}\rightarrow\mathcal{Q}$ must have infinite cardinality.  
\end{remark}

\medskip
\noindent\textbf{Acknowledgements}. We want to thank Gabe Cunningham for 
 insightful comments. We are also particularly grateful to 
the Instituto de Matem\'{a}ticas de la UNAM and allied sponsors
 for hosting the \textit{Second Workshop  on Abstract
Polytopes},  Cuernavaca, Mexico, July 30-August 3, 2012. 
 
\bibliography{MonsSch}

\end{document}